\documentclass[reqno,11pt]{amsart}
\usepackage{amssymb,enumerate}
\usepackage{amsmath}
\usepackage{mathrsfs}
\usepackage{dsfont} 
\usepackage{hyperref}
\usepackage{dsfont}
\usepackage[usenames]{color}
\usepackage{graphicx,graphics}

\allowdisplaybreaks

\begin{document}

\renewcommand{\AA}{\mathcal{A}}
\renewcommand{\div}{{\rm div \,}}
\newcommand{\BB}{\mathcal{B}}
\newcommand{\CC}{\mathcal{C}}
\newcommand{\DD}{\mathcal{D}}
\newcommand{\EE}{\mathcal{E}}
\newcommand{\FF}{\mathcal{F}}
\newcommand{\GG}{\mathcal{G}}
\newcommand{\HH}{\mathcal{H}}
\newcommand{\II}{ I\hspace{-0.5mm}I}
\newcommand{\JJ}{\mathcal{J}}
\newcommand{\KK}{\mathcal{K}}
\newcommand{\LL}{\mathcal{L}}
\newcommand{\MM}{\mathcal{M}}
\newcommand{\NN}{\mathcal{N}}
\newcommand{\OO}{\mathcal{O}}
\newcommand{\PP}{\mathcal{P}}
\newcommand{\QQ}{\mathcal{Q}}
\newcommand{\RR}{\mathcal{R}}
\renewcommand{\SS}{\mathcal{S}}
\newcommand{\TT}{\mathcal{T}}
\newcommand{\UU}{\mathcal{U}}
\newcommand{\VV}{\mathcal{V}}
\newcommand{\WW}{\mathcal{W}}
\newcommand{\XX}{\mathcal{X}}
\newcommand{\YY}{\mathcal{Y}}
\newcommand{\ZZ}{\mathcal{Z}}


\newcommand{\ds}{\displaystyle}

\newcommand{\A}{\mathbb{A}}
\newcommand{\B}{\mathbb{B}}
\newcommand{\C}{\mathbb{C}}
\newcommand{\D}{\mathbb{D}}
\newcommand{\E}{\mathbb{E}}
\newcommand{\F}{\mathbb{F}}
\newcommand{\G}{\mathbb{G}}
\renewcommand{\H}{\mathbb{H}}
\newcommand{\I}{\mathbb{I}}
\newcommand{\J}{\mathbb{J}}
\newcommand{\K}{\mathbb{K}}
\renewcommand{\L}{\mathbb{L}}
\newcommand{\M}{\mathbb{M}}
\newcommand{\N}{\mathbb{N}}
\renewcommand{\O}{\mathbb{O}}
\renewcommand{\P}{\mathbb{P}}
\newcommand{\Q}{\mathbb{Q}}
\newcommand{\R}{\mathbb{R}}
\renewcommand{\S}{\mathbb{S}}
\newcommand{\T}{\mathbb{T}}
\newcommand{\U}{\mathbb{U}}
\newcommand{\V}{\mathbb{V}}
\newcommand{\W}{\mathbb{W}}
\newcommand{\X}{\mathbb{X}}
\newcommand{\Y}{\mathbb{Y}}
\newcommand{\Z}{\mathbb{Z}}

\newcommand{\cC}{{\mathcal C}}


\newcommand{\al}{\alpha}
\newcommand{\be}{\beta}
\newcommand{\ga}{\gamma}
\newcommand{\de}{\delta}
\newcommand{\ep}{\varepsilon}
\newcommand{\ze}{\zeta}
\newcommand{\et}{\eta}
\newcommand{\vth}{\vartheta}
\renewcommand{\th}{\theta}
\newcommand{\io}{\iota}
\newcommand{\ka}{\kappa}
\newcommand{\la}{\lambda}
\newcommand{\rh}{\rho}
\newcommand{\si}{\sigma}
\newcommand{\ta}{\tau}
\newcommand{\up}{\upsilon}
\newcommand{\ph}{\varphi}
\newcommand{\ch}{\chi}
\newcommand{\ps}{\psi}
\newcommand{\om}{\omega}

\newcommand{\Ga}{\Gamma}
\newcommand{\De}{\Delta}
\newcommand{\Th}{\Theta}
\newcommand{\La}{\Lambda}
\newcommand{\Si}{\Sigma}
\newcommand{\Up}{\Upsilon}
\newcommand{\Ph}{\xi}
\newcommand{\Om}{\Omega}

\newcommand{\BR}{\color{red}}
\newcommand{\ER}{\color{black}}

\newcommand{\inj}{\hookrightarrow}
\newcommand{\stetein}{\overset{s}{\hookrightarrow}}
\newcommand{\dichein}{\overset{d}{\hookrightarrow}}
\newcommand{\pa}{\partial}
\newcommand{\re}{\restriction}
\newcommand{\tief}{\downharpoonright}

\newcommand{\bra}{\langle}
\newcommand{\ket}{\rangle}
\newcommand{\bs}{\backslash}
\newcommand{\divv}{\operatorname{div}}
\newcommand{\Dt}{\frac{\mathrm d}{\mathrm dt}}

\newcommand{\sm}{\setminus}
\newcommand{\es}{\emptyset}

\newtheorem{theorem}{Theorem}
\newtheorem{corollary}{Corollary}
\newtheorem*{main}{Main Theorem}
\newtheorem{lemma}[theorem]{Lemma}
\newtheorem{proposition}{Proposition}
\newtheorem{conjecture}{Conjecture}
\newtheorem*{problem}{Problem}
\theoremstyle{definition}
\newtheorem{definition}[theorem]{Definition}
\newtheorem{remark}{Remark}
\newtheorem*{notation}{Notation}

\newcommand{\cqfd}{\begin{flushright}\vspace*{-3mm}$\Box $\vspace{-2mm}\end{flushright}}
\newcommand{\saut}{\vspace*{1mm}\\ \nd }
\newcommand{\on}{\mbox{ on }}
\newcommand{\with}{\mbox{ with }}
\newcommand{\nd}{\noindent}
\newcommand{\eps}{\varepsilon}
\newcommand{\limit}{\underset{n\rightarrow +\infty}{\longrightarrow}}

\sloppy
\title[On the Interface Formation Model]
{On the Interface Formation Model for Dynamic Triple Lines}
\author[D.\ Bothe]{Dieter Bothe}
\address{Center of Smart Interfaces \& Fachbereich Mathematik, Technische Universit\"at Darmstadt, Alarich-Weiss-Str.\ 10, 64287 Darmstadt, Germany}
\email{bothe@csi.tu-darmstadt.de}

\author[J.\ Pr\"uss]{Jan Pr\"uss}
\address{Institut f\"ur Mathematik, Universit\"at Halle-Wittenberg, Theodor-Lieser-Str.\ 5,
06120 Halle (Saale), Germany}
\email{jan.pruess@mathematik.uni-halle.de}

\date{\today}

\begin{abstract}
This paper revisits the theory of Y.\ Shikhmurzaev on forming interfaces as a continuum thermodynamical model for dynamic
triple lines. We start with the derivation of the balances for mass, momentum, energy and entropy in a three-phase fluid
system with full interfacial physics, including a brief review of the relevant transport theorems on interfaces and
triple lines. Employing the entropy principle in the form given in \cite{BotheDreyer}, but extended to this more general case,
we arrive at the entropy production and perform a linear closure, except for a nonlinear closure for the sorption processes.
Specialized to the isothermal case, we obtain a thermodynamically consistent mathematical model for dynamic triple lines
and show that the total available energy is a strict Lyapunov function for this system.
\end{abstract}
\maketitle

\section{Introduction}
The line at which three phases meet, is called a triple line; cf.\ Figure~1.
If the phases which touch each other are all fluid phases, i.e.\
two immiscible liquids are in contact with another liquid or a gas, this triple line is freely deformable in space, while it is
bound to move on a given surface, if one of the phases is a solid. In the latter case, one usually speaks about a dynamic
contact line, while the notion of a triple line is typically used in the former setting. Both cases share many similarities and
their modeling and analysis is closely related. In applications, wetting more often appears on a solid wall, i.e.\
the case of a contact line is more often considered. Hence, the main body of the literature is devoted to this case.
The present paper deals with dynamics triple lines, but in such a generality, that analogous results are valid for the
contact line situation. Nevertheless, due to the more frequent encounter of wetting of solid supports, the brief
literature survey to follow necessarily focuses on contact line dynamics.

The modeling and computation of dynamic contact lines is an active field due to the enormous relevance of wetting and dewetting
phenomena in various technical and industrial applications; see \cite{blake2006physics}, \cite{bonn2009wetting} and
\cite{Snoeijer-CL} for recent surveys on the field, containing also references to experimental work.
Different modeling approaches are employed, containing in particular so-called molecular-kinetic theory (MKT;
see, e.g., \cite{blake1969kinetics}, \cite{de1985wetting}) and continuum physical theories. The latter is often subsumed under the heading ``hydrodynamic theory'' and is mostly based on sharp-interface models, while phase field models have also been extended to
cover contact lines as in \cite{Jacqmin-DI}. The sharp-interface hydrodynamic theory started essentially with the
seminal paper by Huh and Scriven \cite{huh1971hydrodynamic} in which the fundamental problem of the inconsistency between
a moving contact line and a no-slip condition at the fixed wall has been analyzed and shown to lead to a non-integrable
stress singularity; cf.\ also \cite{dussan1974motion} and, for a more rigorous mathematical treatment,
\cite{Pukhnachev83} and \cite{solonnikov95}.
Consequently, subsequent models always rely on some ``relaxation'' at the contact line, and the most common
way to remove the stress singularity, as already proposed in \cite{huh1971hydrodynamic}, is to introduce Navier-type slip close to the contact line. Besides this complication, the main extension of the standard two-phase Navier-Stokes system consists of
the prescription of the dynamic contact angle $\theta_d$, i.e.\ the angle which is formed between the fluid interface and the solid support,
as a function of the contact line speed. At this point it is to be noted that the contact angle changes its value under
dynamic conditions, while the equilibrium contact angle $\theta_e$ is usually assumed to be governed by Young's law, i.e.\
\begin{equation}
\sigma^{gl} \cos (\theta_e) = \sigma^{gs} -\sigma^{ls}
\end{equation}
in case of a liquid wetting a solid surrounded by a gas phase,
where the superscripts stand for gas (g), liquid (l) and solid (s) and $\sigma$ denotes the interfacial tension
of the respective interface. Based on the classical
experimental studies in \cite{hoffman1975study}, the general form of the relation between $\theta_d$ and the contact line speed
is given by the heuristic relation
\begin{equation}
\theta_d = f_{\rm Hoff} \big( {\rm Ca} + f_{\rm Hoff}^{-1} (\theta_e )\big),
\end{equation}
where Ca denotes the Capillary number given as ${\rm Ca}= \eta U / \sigma^{gl}$ with $\eta$ the dynamic liquid
viscosity and $U$ the contact line speed. Several concrete correlations have been established for different materials and
certain wetting scenarios like ``Tanner's law'' \cite{tanner1979spreading} or the correlation of Jiang et al.\ \cite{jiang1979correlation}. Theoretical investigations using the hydrodynamic theory identified three length-scales near the
contact line: an inner region in which the fluid interface is essentially planar and touches the solid support at the equilibrium
angle; a mesoscopic region in which a significant bending of the interface can occur; an outer (macroscopic) region in which the
contact angle attains a different value, the so-called apparent contact angle.
The hydrodynamic theory provides relations for the dependence
of the contact angle on the distance from the contact line especially in the mesoscopic region; see \cite{voinovhydrodynamics}, \cite{coxdynamics}, \cite{dussan1976moving}.
Knowledge of this dependence is very useful for numerical purpose, both as a subgrid-scale model to reduce the necessary
resolution at the contact-line and in order to
neutralize the inherent mesh dependence of numerical solutions due to the typical under-resolution of the smallest length scales
in the contact line region; for the latter, see \cite{Afkhami} and \cite{fath-CL}.

\begin{figure}[t]
\centering
\includegraphics[width=4.6in]{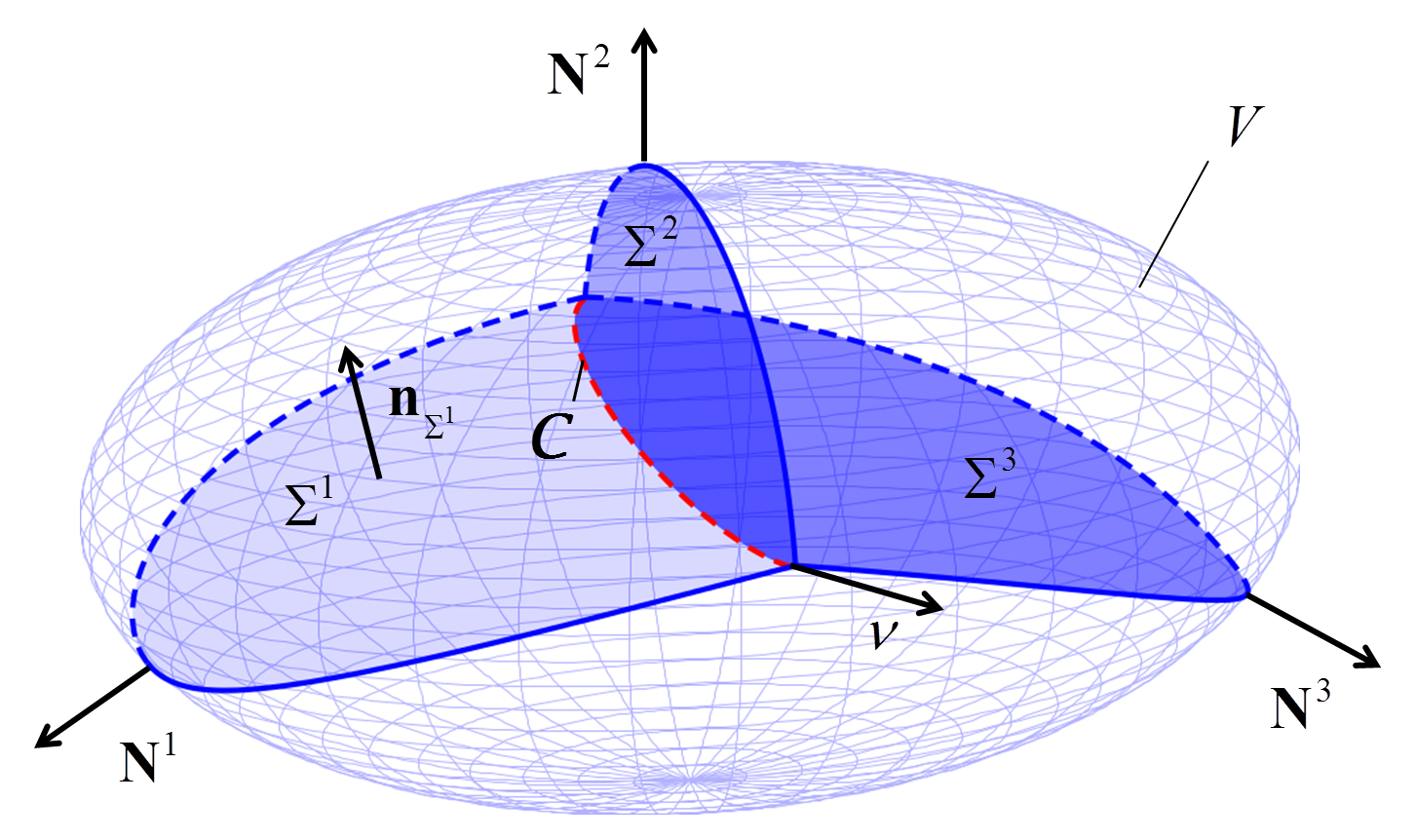}
\caption{Configuration of the phases and interfaces at the contact line.}
\end{figure}
While the hydrodynamic model can describe many wetting processes at least qualitatively, in particular concerning the observed
dynamical shapes of attached droplets moving on a wall, say, it does not capture the full physics of a dynamic contact line.
One important deviation is the internal flow field in the wetting liquid close to an advancing contact line, which experimentally
is known to be a rolling motion (\cite{dussan1974motion}, \cite{chen1997velocity}), but is a sliding motion in numerical
simulations using the above model. Moreover, there is experimental evidence that the relation between the dynamic contact
angle and the contact line speed is more complicated and not of such a simple local nature; cf.\ \cite{blake1999experimental}.
For further discrepancies between experimental observations and the hydrodynamic model see \cite{Shikh-book}.
A very interesting approach to overcome these short-comings has been introduced by Y.\ Shikhmurzaev in \cite{Shikh93}; see
also \cite{Shikh-book}.
The approach there also employs continuum physics, but accounts for the aspect of interface formation and disappearance
at the contact line. A crucial point for the model development then is to include enough interfacial thermodynamics to
allow for a non-constant interfacial tension, governed by a surface equations of state on all involved interfaces.
For this purpose, the mass contained in the interfacial layer has to be balanced separately, since it encounters different forces
compared to within the bulk phases and it is this mass density which determines the surface pressure, i.e.\ the surface tension.
In the considered sharp interface/sharp contact line model, the interfacial mass is lumped into an area-specific mass
density and the model is extended to cover the evolution of this interfacial mass density by appropriate balance equations
on the moving surfaces. This model has proven a great potential to explain several physical phenomena like wetting, coalescence, cusp formation and the break-up of liquid threads; cf.\ \cite{Shikh-book}, \cite{blake2006physics} and the introduction in
\cite{billingham06}.
The interface formation model of Shikhmurzaev
has been based on the continuum thermodynamics of fluid interfaces developed in \cite{BAM},
\cite{Bedeaux}, but with the sensible aim to formulate the most simple model which is able to describe the wetting process with dynamic contact angle and rolling motion close to the contact line with the material properties modeled via bulk and surface free energies but without a heuristic relation between contact angle and contact line speed.
Several years after the fundamental paper \cite{Shikh93} appeared, Billingham in \cite{billingham06} pointed out that one
further condition at the contact line has to be added, and he employed a condition provided by Bedeaux in \cite{Bedeaux-CL}
which relates the rate of mass transfer from one surface into the other to the difference of the surface chemical potentials.
We will come back to this point in the final remarks at the end of this paper.

In the brief survey above, the topic called ``contact angle hysteresis'', referring to the appearance of a full interval
of possible contact angles in the static case which spans the range from the angles observed for (infinitely slowly)
advancing and receding contact lines, has not been touched.
This phenomenon seems to be similar to dry friction between solids and, in fact,
the notion of contact line friction is also present in the literature on the molecular kinetic theory of contact lines.
For this topic, we refer to \cite{bonn2009wetting}, \cite{Snoeijer-CL} and the references given there.

\section{Integral Balances}
We consider a region $G \subset \R^3$ filled with three bulk phases $\Omega^k (t)$ ($k=1, 2, 3$), separated by interfaces
$\Sigma^k (t)$ ($k=1, 2, 3$) which meet at a common triple line $\mathcal C$. As an example, imagine a liquid phase $\Omega^1$ in
the form of a water droplet sitting on another liquid, say oil, which forms bulk phase $\Omega^2$, and being surrounded by phase
$\Omega^3$ composed of air. Then, for instance, $\Sigma^1$ denotes the oil-water interface, $\Sigma^2$ the interface
between the air and the oil and $\Sigma^3$ the air-water interface. The deformable and free
bounding curve at which all three interfaces meet is the so-called triple line ${\mathcal C}$.
As a related but somewhat different case, consider again a liquid phase $\Omega^1$ in
the form of a droplet, but now sitting on a solid support, which forms bulk phase $\Omega^2$, and being surrounded by gas phase
$\Omega^3$. Then two out of the three interfaces are fixed and the triple line is the set of all points
where the gas-liquid interface meets the solid support. In this case, one usually calls ${\mathcal C}$ the contact line
which now has reduced degrees of freedom due to the solid support.
We focus on three-phase fluid systems with a common triple line and assume that the interfaces meet at angles different from
0 and $\pi$; we shall refer to this as the non-degenerate case.

We start with the integral balance of a generic extensive quantity which is present in the bulk phases with specific density $\phi$, on the interfaces with specific density $\phi^\Sigma$ and on the triple line with specific density $\phi^{\mathcal C}$.
Hence $\rho \phi$, $\rho^\Sigma \phi^\Sigma$ and $\rho^{\mathcal C} \phi^{\mathcal C}$, respectively, are the
volume-, area- and line-specific densities, where $\rho $, $\rho^\Sigma $ and $\rho^{\mathcal C}$ are the mass densities.
If a specific bulk phase or interface is considered, we write $\rho_k \phi_k$ or $\rho^\Sigma_k \phi^\Sigma_k$, respectively, for the respective density.
With this notation, the generic integral balance for a fixed control volume $V\subset G$ reads as
\begin{align}
& \frac{d}{dt} \left[ \int_{\Omega_V} \rho \phi \,dx + \int_{\Sigma_V} \rho^\Sigma \phi^\Sigma \,do
+ \int_{\mathcal{C}_V} \rho^{\mathcal C} \phi^{\mathcal C} \,dl \right] = \nonumber \\[1ex]
& - \int_{\partial \Omega_V} (\rho \phi v + j) \cdot n\,do
- \int_{\partial \Sigma_V} (\rho^\Sigma \phi^\Sigma v^\Sigma + j^\Sigma) \cdot N\,dl
- \int_{\partial \mathcal{C}_V} (\rho^{\mathcal C} \phi^{\mathcal C} v^{\mathcal C} + j^{\mathcal C})
\cdot \nu \,dP \label{generic-balance}\\[1ex]
& + \int_{\Omega_V} f \,dx + \int_{\Sigma_V} f^\Sigma \,do + \int_{\mathcal{C}_V} f^{\mathcal C} \,dl .\nonumber
\end{align}
In \eqref{generic-balance} we let $dx,\, do$ and $dl$ denote the volume, area and line measure, respectively. Moreover,
$dP$ denotes the point (i.e., counting) measure.

Here as well as throughout the paper we use the following condensed notation. First,
\[
\Omega_V := \bigcup_{k=1}^3 \Omega_{V}^k, \quad
\Sigma_V := \bigcup_{k=1}^3 \Sigma_{V}^k
\quad \mbox{ with } \Omega_V^k:=\Omega^k \cap V, \quad
\Sigma_V^k:=\Sigma^{k}\cap V,
\]
which are all time-dependent sets. We assume a single triple line, hence $\mathcal{C}_V:=\mathcal{C} \cap V$. The detailed
version of \eqref{generic-balance} then reads as
\begin{align*}
& \frac{d}{dt} \left[ \sum_{k=1}^3 \int_{\Omega_V^k} \rho_k \phi_k \,dx \; + \;
\sum_{k=1}^3 \int_{\Sigma_V^k} \rho^\Sigma_k \phi^\Sigma_k \,do
+ \int_{\mathcal{C}_V} \rho^{\mathcal C} \phi^{\mathcal C} \,dl \right] =  \\[1ex]
& - \sum_{k=1}^3  \int_{\partial \Omega_V^k} (\rho_k \phi_k v_k + j_k) \cdot n^k\,do
- \sum_{k=1}^3 \int_{\partial \Sigma_V^k} (\rho^\Sigma_k \phi^\Sigma_k v^\Sigma_k + j^\Sigma_k) \cdot N^k\,dl
- \int_{\partial \mathcal{C}_V} (\rho^{\mathcal C} \phi^{\mathcal C} v^{\mathcal C} + j^{\mathcal C})
\cdot \nu \,dP  \\[1ex]
& + \sum_{k=1}^3 \int_{\Omega_V^k} f_k \,dx + \sum_{k=1}^3 \int_{\Sigma_V^k} f^\Sigma_k \,do + \int_{\mathcal{C}_V} f^{\mathcal C} \,dl .
\end{align*}
For better readability, we use the condensed notation whenever this is reasonable.

We apply this balancing to the extensive quantities mass, momentum, energy and entropy. The corresponding integral balances read as\\[1ex]
{\bf mass balance.}\\
\begin{align}
 \frac{d}{dt} \left[ \int_{\Omega_V} \rho \,dx + \int_{\Sigma_V} \rho^\Sigma \,do
+ \int_{\mathcal{C}_V} \rho^{\mathcal C} \,dl \right] = 
 - \int_{\partial \Omega_V} \rho v \cdot n\,do
- \int_{\partial \Sigma_V} \rho^\Sigma  v^\Sigma \cdot N\,dl
- \int_{\partial \mathcal{C}_V} \rho^{\mathcal C} v^{\mathcal C} \cdot \nu \,dP.
\end{align}
$\mbox{ }$\\[1ex]
{\bf momentum balance.}\\
\begin{align}
& \frac{d}{dt} \left[ \int_{\Omega_V} \rho v \,dx + \int_{\Sigma_V} \rho^\Sigma v^\Sigma \,do
+ \int_{\mathcal{C}_V} \rho^{\mathcal C} v^{\mathcal C} \,dl \right] = \nonumber \\[1ex]
& - \int_{\partial \Omega_V} \rho v (v \cdot n)\,do
- \int_{\partial \Sigma_V} \rho^\Sigma  v^\Sigma (v^\Sigma\cdot N)\,dl
- \int_{\partial \mathcal{C}_V} \rho^{\mathcal C} v^{\mathcal C} (v^{\mathcal C}\cdot \nu) \,dP
 \\[1ex]
& + \int_{\partial \Omega_V} S \cdot n\,do
+ \int_{\partial \Sigma_V} S^\Sigma \cdot N\,dl
+ \int_{\partial \mathcal{C}_V} S^{\mathcal C}\cdot \nu \,dP
 + \int_{\Omega_V} \rho b \,dx + \int_{\Sigma_V} \rho^\Sigma b^\Sigma \,do
 + \int_{\mathcal{C}_V} \rho^{\mathcal C} b^{\mathcal C} \,dl .\nonumber
\end{align}
Here $S$, $S^\Sigma $ and $S^{\mathcal C}$ denote the stress tensor in the bulk phases,
on the interfaces and on the triple line, respectively, and $b$, $b^\Sigma $ and $b^{\mathcal C}$
are the specific body forces.\\[2ex]
{\bf energy balance.}\\
\begin{align}
\label{total-energy}
& \frac{d}{dt} \left[ \int_{\Omega_V} \rho (e+\frac{v^2}{2}) \,dx
+ \int_{\Sigma_V} \rho^\Sigma (e^\Sigma+ \frac{(v^\Sigma)^2}{2}) \,do
+ \int_{\mathcal{C}_V} \rho^{\mathcal C} (e^{\mathcal C} + \frac{(v^{\mathcal C})^2}{2}) \,dl \right] = \nonumber \\[1ex]
& - \int_{\partial \Omega_V} \rho (e+\frac{v^2}{2}) v \cdot n\,do
- \int_{\partial \Sigma_V} \rho^\Sigma  (e^\Sigma+ \frac{(v^\Sigma)^2}{2}) v^\Sigma\cdot N\,dl
- \int_{\partial \mathcal{C}_V} \rho^{\mathcal C} (e^{\mathcal C} + \frac{(v^{\mathcal C})^2}{2}) v^{\mathcal C}\cdot \nu \,dP
 \\[1ex]
& + \int_{\partial \Omega_V} (v\cdot S -q) \cdot n\,do
+ \int_{\partial \Sigma_V} (v^\Sigma \cdot S^\Sigma - q^\Sigma) \cdot N\,dl
+ \int_{\partial \mathcal{C}_V} (v^{\mathcal C} \cdot S^{\mathcal C}- q^{\mathcal C}) \cdot \nu \,dP\nonumber \\[1ex]
& + \int_{\Omega_V} \rho v\cdot b \,dx + \int_{\Sigma_V} \rho^\Sigma v^\Sigma \cdot b^\Sigma \,do
 + \int_{\mathcal{C}_V} \rho^{\mathcal C} v^{\mathcal C} \cdot b^{\mathcal C} \,dl .\nonumber
\end{align}
Here $q$, $q^\Sigma $ and $q^{\mathcal C}$ denote the heat flux in the bulk phases,
on the interfaces and on the triple line, respectively. Note that energy sources due to radiation have been
omitted in \eqref{total-energy}.\\[2ex]
{\bf entropy balance.}\\
\begin{align}
\label{entropy-integral}
& \frac{d}{dt} \left[ \int_{\Omega_V} \rho s \,dx
+ \int_{\Sigma_V} \rho^\Sigma s^\Sigma \,do
+ \int_{\mathcal{C}_V} \rho^{\mathcal C} s^{\mathcal C} \,dl \right] = \nonumber \\[1ex]
& - \int_{\partial \Omega_V} (\rho s v +\Phi) \cdot n\,do
- \int_{\partial \Sigma_V} (\rho^\Sigma  s^\Sigma  v^\Sigma +\Phi^\Sigma) \cdot N\,dl
- \int_{\partial \mathcal{C}_V} (\rho^{\mathcal C} s^{\mathcal C} v^{\mathcal C} + \Phi^{\mathcal C})\cdot \nu \,dP
 \\[1ex]
& + \int_{\Omega_V} \zeta \,dx + \int_{\Sigma_V} \zeta^\Sigma \,do
 + \int_{\mathcal{C}_V} \zeta^{\mathcal C} \,dl .\nonumber
\end{align}
Here $\Phi$, $\Phi^\Sigma $ and $\Phi^{\mathcal C}$ denote the entropy flux in the bulk phases,
on the interfaces and on the triple line, respectively, while $\zeta$, $\zeta^\Sigma $ and $\zeta^{\mathcal C}$
are the corresponding entropy productions. \\[1ex]
\noindent
{\bf Remark.} Note that the internal energy density as well as the entropy density can be positive even if the area-
or line-specific mass densities are considered to be zero. In other words, in the limit as $\rho^\Sigma \to 0+$
or $\rho^{\mathcal C} \to 0+$, products such as $\rho^\Sigma e^\Sigma$ or $\rho^{\mathcal C} s^{\mathcal C}$ may converge
to strictly positive limit densities, i.e.\
\[
\rho^\Sigma e^\Sigma \to u^\Sigma, \quad
\rho^\Sigma s^\Sigma \to \eta^\Sigma, \quad
\rho^{\mathcal C} e^{\mathcal C} \to u^{\mathcal C}, \quad
\rho^{\mathcal C} s^{\mathcal C} \to \eta^{\mathcal C}
\]
with non-vanishing densities $u^\Sigma, \eta^\Sigma, u^{\mathcal C}, \eta^{\mathcal C}$ has to be allowed for.
Otherwise, for instance, the surface tension for a fluid interface with zero surface mass density would automatically vanish.

\section{Transport Theorems}
The derivation of local versions of the balance equations follow by application of appropriate transport theorems and subsequent localization. The following transport theorems will be employed.\\[1ex]
{\bf Volume transport.}
In the general setting described above, let $V\subset \R^3$ be a fixed control volume in $G$, let $\Sigma$ be short for $\bigcup_{k=1}^3 \Sigma^k$ with the time-dependent interfaces $\Sigma^k (t)$ and $n_\Sigma=n_{\Sigma^k}$ the unit normal field on $\Sigma^k(t)$ with an
arbitrary fixed orientation. Let $V_\Sigma$ denote the speed of normal displacement of $\Sigma^k( \cdot)$. The latter is
a purely kinematic quantity, but it is related to the barycentric velocity of the interfacial mass via
$V_\Sigma= v^\Sigma \cdot n_\Sigma$. Moreover, given any bulk field $\phi$, the jump of $\phi$ at $\Sigma$ is defined
by the jump bracket $[\![\cdot ]\!]$ according to
\begin{equation}
[\![\phi ]\!](t,x):= \lim_{h\to 0+} \big( \phi(t,x+h n_\Sigma )-  \phi(t,x - h n_\Sigma ) \big).
\end{equation}
With these notations and for the specific control volumes mentioned as well as for sufficiently smooth fields, it holds that
\begin{equation}
\label{vol-trans}
\frac{d}{dt} \int_{V} \phi \,dx = \int_{V\setminus \Sigma} \partial_t \phi \,dx
- \int_{\Sigma_V} [\![ \phi ]\!] \, V_\Sigma \,do,
\end{equation}
where $\Sigma_V (t):=\Sigma(t) \cap V$.\\[1ex]
{\bf Surface transport.}
In the general setting described above, let $V\subset \R^3$ be a fixed control volume in $G$. Then, for sufficiently smooth fields, it holds that
\begin{equation}
\label{surf-trans}
\frac{d}{dt} \int_{\Sigma_V} \phi^\Sigma \,do
= \int_{\Sigma_V} \big( \partial_t^\Sigma \phi^\Sigma - \phi^\Sigma \kappa_\Sigma V_\Sigma \big) \,do
+ \int_{\partial \Sigma_V} \phi^\Sigma \, V_{\partial \Sigma_V} \,dl.
\end{equation}
Here $\partial_t^\Sigma$ denotes the time derivative along a path that follows the normal motion of $\Sigma (\cdot )$ and
$\kappa_\Sigma := \div_\Sigma (- n_\Sigma)$ is twice the mean curvature. Furthermore,
$V_{\partial \Sigma_V}$ is the normal (relative to the boundary of $\Sigma_V$) speed of displacement of $\partial \Sigma_V (\cdot)$ (in the plane tangential to $\Sigma$).

Let us note in passing that the derivation of the local balance equations can be done with special control volumes
such that the outer normal $n_V$ satisfies $n_V \perp n_\Sigma$ on $\partial V \cap \Sigma$. For such control volumes
the boundary contribution, i.e.\ the final term in \eqref{surf-trans}, vanishes.\\[1ex]
{\bf Line transport.}
For sufficiently smooth fields, it holds that
\begin{equation}
\label{CL-trans}
\frac{d}{dt} \int_{{\mathcal C}_V} \phi^{\mathcal C} \,dl
= \int_{{\mathcal C}_V} \big( \frac{D^{\mathcal C} \phi^{\mathcal C}}{Dt} + \phi^{\mathcal C}
\div_{\mathcal C} v^{\mathcal C} \big) \,dl
+ \int_{\partial {\mathcal C}_V} \phi^{\mathcal C} \, \big(V_{\partial {\mathcal C}_V} - v^{\mathcal C} \cdot \nu \big)\,dP.
\end{equation}
Here $\frac{D^{\mathcal C}}{Dt}$ denotes the Lagrangian derivative, following the triple line along a path with velocity
$v^{\mathcal C}$ and $V_{\partial {\mathcal C}_V}$ is the normal (relative to the end points of ${\mathcal C}_V$) speed
of displacement of $\partial {\mathcal C}_V$. Recall that $\nu$ is the outer normal to the curve ${\mathcal C}_V$ in its
end points (cf.\ Figure 1) and that $dP$ denotes the point (i.e., counting) measure.\\[2ex]
{\bf Remarks.}
1. The transport theorems above appear rather different.
Actually, they can all be brought into the same form as the line transport theorem.
In case of surface transport, this follows directly from the relation
\[
\frac{D^\Sigma \phi^\Sigma}{Dt} = \partial_t^\Sigma \phi^\Sigma + v^\Sigma \cdot \nabla_\Sigma \phi^\Sigma
\]
for the surface Lagrangian derivative. Note that $\div_\Sigma v^\Sigma = \div_\Sigma v^\Sigma_{||} - \kappa_\Sigma V_\Sigma$,
hence
\[
\partial_t^\Sigma \phi^\Sigma - \phi^\Sigma \kappa_\Sigma V_\Sigma=
\frac{D^\Sigma \phi^\Sigma}{Dt} - \div_\Sigma (\phi^\Sigma v^\Sigma_{||}) + \phi^\Sigma \div_\Sigma v^\Sigma
\]
and then, by the surface divergence theorem, equation \eqref{surf-trans} implies
\begin{equation}
\label{surf-trans-general}
\frac{d}{dt} \int_{\Sigma_V} \phi^\Sigma \,do
= \int_{\Sigma_V} \big( \frac{D^\Sigma \phi^\Sigma}{Dt} + \phi^\Sigma \div_\Sigma v^\Sigma \big) \,do
+ \int_{\partial \Sigma_V} \phi^\Sigma \, \big( V_{\partial \Sigma_V}- v^\Sigma \cdot N \big) \,dl.
\end{equation}
To bring the volume transport formula \eqref{vol-trans} into the same form, one first observes that \eqref{vol-trans}
combines the transport formulas for both bulk phases which meet at the considered interface.
If two bulk phases $\Omega^\pm (t)$ are separated by an interface $\Sigma (t)$, then a simple variant of the Reynolds
transport theorem yields
\begin{equation}
\frac{d}{dt} \int_{\Omega^\pm_V} \phi  \,dx
= \int_{\Omega^\pm_V} \big( \frac{D \phi}{Dt} + \phi \, \div v  \big) \,dx
+ \int_{\partial \Omega^\pm_V} \phi \, \big( V_{\partial \Omega^\pm_V}- v  \cdot n \big) \,do,
\end{equation}
where $n$ is the outer unit normal to $\Omega^\pm_V$; note that the latter coincides with $\pm n_\Sigma$
on $\Sigma_V$.\\[1ex]
\indent
2. An equivalent form of \eqref{surf-trans-general} reads as
\begin{equation}
\label{surf-trans-2}
\frac{d}{dt} \int_{\Sigma_V} \phi^\Sigma \,do
= \int_{\Sigma_V} \big( \frac{D^\Sigma \phi^\Sigma}{Dt} + \phi^\Sigma \div_\Sigma v^\Sigma \big) \,do
- \int_{\partial \Sigma_V} \phi^\Sigma \, \frac{v^\Sigma \cdot n_V}{\sqrt{1-(n_{\Sigma} \cdot n_V)^2}} \,dl,
\end{equation}
where $n_V$ is the outer unit normal to $V$. For this purpose, one first uses elementary geometry to
compute $V_{\partial \Sigma_V}=-V_\Sigma \frac{n_\Sigma \cdot n_V}{\sqrt{1-(n_{\Sigma} \cdot n_V)^2}}$.
Since $\sqrt{1-(n_{\Sigma} \cdot n_V)^2}=N\cdot n_V$ and $\{N,n_\Sigma, \tau \}$ with $\tau$ a unit vector tangential
to $\partial \Sigma_V$ (hence also to $\partial V$) is a local orthonormal basis, the equation \eqref{surf-trans-2} follows from
\[
 v^\Sigma \cdot N - V_{\partial \Sigma_V} = \frac{1}{\sqrt{1-(n_{\Sigma} \cdot n_V)^2}}
 \Big( (v^\Sigma \cdot N)\, (N\cdot n_V) + (v^\Sigma \cdot n_\Sigma)\, (n_\Sigma\cdot n_V) \Big).
\]
The relation from \eqref{surf-trans-general} has been given in \cite{CFG}, while the variant \eqref{surf-trans-2}
can be found in Chapter~3 in \cite{Romano}; see also the appendix in \cite{AB}.\\[1ex]
\indent
3. Below we will also use variants of the above transport theorems with built-in mass balance. These read as
\begin{equation}
\label{vol-trans2}
\frac{d}{dt} \int_{V} \rho \phi \,dx = \int_{V\setminus \Sigma} \rho \frac{D \phi}{Dt} \,dx
+ \int_{\Sigma_V} [\![ \dot{m} \phi ]\!] \,do - \int_{\partial V} \rho \phi v \cdot n \,do
\end{equation}
with $\dot{m}^\pm:=\rho^\pm (v^\pm-v^\Sigma )\cdot n_\Sigma$ on $\Sigma$, and
\begin{equation}
\label{surf-trans2}
\frac{d}{dt} \int_{\Sigma_V} \rho^\Sigma \phi^\Sigma \,do
= \int_{\Sigma_V} \big( \rho^\Sigma \frac{D^\Sigma \phi^\Sigma}{Dt} - [\![ \dot{m} \phi^\Sigma ]\!] \big) \,do
+ \int_{\partial \Sigma_V} \rho^\Sigma \phi^\Sigma \, \big( V_{\partial \Sigma_V}- v^\Sigma \cdot N \big) \,dl.
\end{equation}

\section{Local Balances}
Application of the transport theorems and localization yields the following local balance equations\\[1ex]
{\bf bulk phase.}
\begin{align}
& \partial_t \rho + \div (\rho v) =0,  \label{bulk-mass}\\[1ex]
& \partial_t (\rho v) + \div (\rho v\otimes v-S)=\rho b, \label{bulk-momentum}\\[1ex]
& \partial_t (\rho e) + \div (\rho e v + q)=S : \nabla v, \hspace{4in}  \\[1ex]
& \partial_t (\rho s) + \div (\rho s v + \Phi)=\zeta .
\end{align}
These are the well-known balance equations in a bulk phase.\\[1ex]
{\bf interface.}
\begin{align}
& \partial_t^\Sigma \rho^\Sigma + \div_\Sigma (\rho^\Sigma v^\Sigma)
+[\![ \rho(v-v^\Sigma)\cdot n_\Sigma ]\!] =0, \label{surf-mass}\\[1ex]
& \partial_t^\Sigma (\rho^\Sigma v^\Sigma)+ \div_\Sigma (\rho^\Sigma v^\Sigma\otimes v^\Sigma -S^\Sigma)
+[\![ \big( \rho v\otimes(v-v^\Sigma)-S \big) \cdot n_\Sigma]\!] =\rho^\Sigma b^\Sigma, \label{surf-momentum}\\[1ex]
&\partial_t^\Sigma (\rho^\Sigma e^\Sigma)+ \div_\Sigma (\rho^\Sigma e^\Sigma  v^\Sigma + q^\Sigma)
+[\![ \Big( \rho (e+\frac{(v-v^\Sigma)^2}{2})(v-v^\Sigma)- (v-v^\Sigma)\cdot S + q \Big) \cdot n_\Sigma]\!]\hspace{0.6in} \nonumber \\
& \qquad = S^\Sigma : \nabla_\Sigma v^\Sigma,  \\[1ex]
& \partial_t^\Sigma (\rho^\Sigma s^\Sigma)+ \div_\Sigma (\rho^\Sigma s^\Sigma v^\Sigma + \Phi^\Sigma)
+[\![ \big( \rho s (v-v^\Sigma)+\Phi \big) \cdot n_\Sigma ]\!] =\zeta^\Sigma  .
\end{align}
Observe that the jump terms always appear with $n_\Sigma$ as a factor. Therefore, these terms are invariant under
re-orientation of the interfaces. Actually, the notion of a ``jump condition'' for these terms can be misleading.
Note that $[\![ f\cdot n_\Sigma ]\!] = - f^+ \cdot n^+ - f^- \cdot n^-$ if the interface separates two bulk phases
$\Omega^\pm$ with outer unit normals $n^\pm$. Hence, if $f$ denotes a bulk flux, the term $-[\![ f\cdot n_\Sigma ]\!]$
describes the total rate of transfer from the bulk phases to the interfaces due to these fluxes.
For the derivation of closure rates below, the explicit form of this term is to be used since two binary
products are involved.\\[1ex]
{\bf triple line.}
\begin{align}
& \partial_t^{\mathcal C} \rho^{\mathcal C} + \div_{\mathcal C} (\rho^{\mathcal C} v^{\mathcal C})
+[\![\![ \rho^\Sigma (v^\Sigma -v^{\mathcal C})\cdot N ]\!]\!] =0,\label{CL-mass} \\[1ex]
& \partial_t^{\mathcal C} (\rho^{\mathcal C} v^{\mathcal C})+ \div_{\mathcal C} (\rho^{\mathcal C} v^{\mathcal C}\otimes v^{\mathcal C} -S^{\mathcal C})
+[\![\![ \big( \rho^\Sigma v^\Sigma\otimes(v^\Sigma-v^{\mathcal C})-S^\Sigma \big) \cdot N]\!]\!] =\rho^{\mathcal C} b^{\mathcal C},\\[1ex]
&\partial_t^{\mathcal C} (\rho^{\mathcal C} e^{\mathcal C})+ \div_{\mathcal C} (\rho^{\mathcal C} e^{\mathcal C}  v^{\mathcal C} + q^{\mathcal C})
+[\![\![ \Big( \rho^\Sigma (e^\Sigma+\frac{(v^\Sigma-v^{\mathcal C})^2}{2})(v^\Sigma-v^{\mathcal C})- (v^\Sigma-v^{\mathcal C})\cdot S^\Sigma + q^\Sigma \Big) \cdot N]\!]\!]\quad \nonumber \\
& \qquad = S^{\mathcal C} : \nabla_{\mathcal C} v^{\mathcal C},  \\[1ex]
& \partial_t^{\mathcal C} (\rho^{\mathcal C} s^{\mathcal C})+ \div_{\mathcal C} (\rho^{\mathcal C} s^{\mathcal C} v^{\mathcal C} + \Phi^{\mathcal C})
+[\![\![ \big( \rho^\Sigma s^\Sigma (v^\Sigma-v^{\mathcal C})+\Phi^\Sigma \big) \cdot N]\!]\!] =\zeta^{\mathcal C}.
\end{align}
Here the triple bracket $[\![\![ \cdot ]\!]\!]$ is defined exclusively for quantities of the form $f^\Sigma \cdot N$ by means of
\begin{equation}
\label{triple-bracket}
[\![\![ f^\Sigma \cdot N ]\!]\!] = - \sum_{k=1}^3 f^{\Sigma}_k \cdot N^k \quad \mbox{ on } {\mathcal C},
\end{equation}
where the sum runs over all interfaces which meet at the triple line and $f^{\Sigma}_k:=f_{|\Sigma^k}$. Let us briefly explain the appearance of such terms,
e.g., for the mass balance \eqref{CL-mass}.
The transport relation \eqref{surf-trans2} for $\phi^\Sigma \equiv 1$ yields the boundary contribution
of the interfacial mass balance as
\[
\int_{\partial \Sigma_V} \rho^\Sigma \, \big( V_{\partial \Sigma_V}- v^\Sigma \cdot N \big) \,dl=
\sum_{k=1}^3 \int_{\partial \Sigma_V^k} \rho^\Sigma_k \, \big( V_{\partial \Sigma_V^k}- v^\Sigma_k \cdot N^k \big) \,dl.
\]
The boundary of $\Sigma_V^k$ is $(\Sigma^k \cap \partial V) \cup {\mathcal C}_V$, hence
\[
\int_{\partial \Sigma_V} \rho^\Sigma \, \big( V_{\partial \Sigma_V}- v^\Sigma \cdot N \big) \,dl=
\sum_{k=1}^3 \int_{\Sigma^k \cap \partial V} \rho^\Sigma_k \, \big( V_{\partial \Sigma_V^k}- v^\Sigma_k \cdot N^k \big) \,dl+
\int_{{\mathcal C}_V} \sum_{k=1}^3 \rho^\Sigma_k \, \big(v^{\mathcal C} - v^\Sigma_k \big) \cdot N^k \,dl.
\]
Employing the condensed notation, this becomes
\[
\int_{\partial \Sigma_V} \rho^\Sigma \, \big( V_{\partial \Sigma_V}- v^\Sigma \cdot N \big) \,dl=
\int_{\Sigma \cap \partial V} \rho^\Sigma \, \big( V_{\partial \Sigma_V}- v^\Sigma \cdot N \big) \,dl
- \int_{{\mathcal C}_V} [\![\![ \rho^\Sigma \, \big(v^{\mathcal C} - v^\Sigma \big) \cdot N ]\!]\!]  \,dl.
\]
%
%

\section{Entropy Production and Closure Relations}
The entropy principle states that every admissible closure for the entropy flux is such that the remaining
entropy production is a sum, running over all dissipative mechanisms, of binary products. The entropy production is
non-negative for any thermodynamic process, i.e.\ the entropy inequality holds. The system is in equilibrium, if and only
if the entropy production vanishes.
For more information about the employed
entropy principle see \cite{BotheDreyer}. We are going to apply this for bulk, interface and triple line in a fully
analogous manner; the details will only be explained for the bulk case. We consider the simplest class of bulk, interface and
contact line materials for which the entropy density is assumed to be a concave function of temperature and mass density, only.
We hence employ constitutive relations of the form
\begin{equation}
\label{entropy-functions}
\rho s = h(\rho e, \rho),\qquad
\rho^\Sigma s^\Sigma = h^\Sigma ( \rho^\Sigma e^\Sigma, \rho^\Sigma ),\qquad
\rho^{\mathcal C} s^{\mathcal C} = h^{\mathcal C} ( \rho^{\mathcal C} e^{\mathcal C}, \rho^{\mathcal C} )
\end{equation}
with concave functions $h$, $h^\Sigma$ and $h^{\mathcal C}$. We furthermore define the (absolute) temperature and
the chemical potential in the respective phase as
\begin{equation}
\label{temperature-defs}
\frac 1 T = \frac{\partial h}{\partial (\rho e)}, \qquad
\frac {1}{T^\Sigma} = \frac{\partial h^\Sigma}{\partial (\rho^\Sigma e^\Sigma)}, \qquad
\frac {1}{T^{\mathcal C}} = \frac{\partial h^{\mathcal C}}{\partial (\rho^{\mathcal C} e^{\mathcal C})}
\end{equation}
and
\begin{equation}
\label{chempot-defs}
-\frac \mu T = \frac{\partial h}{\partial \rho}, \qquad
-\frac{{\mu}^\Sigma}{T^\Sigma} = \frac{\partial h^\Sigma}{\partial \rho^\Sigma}, \qquad
-\frac{{\mu}^{\mathcal C}}{T^{\mathcal C}} = \frac{\partial h^{\mathcal C}}{\partial \rho^{\mathcal C}}.
\end{equation}
We insert the constitutive relation \eqref{entropy-functions} for the entropy density into the respective entropy balance,
use the chain rule employing the definitions \eqref{temperature-defs} and \eqref{chempot-defs} and eliminate all partial
time derivatives by means of the other balance equations. The resulting terms are grouped in such a way that only a single
full divergence appears, which contains in particular the entropy flux, all terms with the velocity divergence as a factor
are collected and all remaining terms are grouped to form a sum of binary products.\\[1ex]
{\bf bulk phase.} The procedure above yields
\begin{align}
\zeta = \div (\Phi - \frac q T )
- \frac{1}{T}(\rho e + P -\rho s T -  \rho \mu)\, \div v
+ q \cdot \nabla \frac{1}{T}
 + \frac{1}{T} S^{\circ}:\nabla v,
\end{align}
where $P:=- \frac 1 3 {\rm tr}\, S$ is the mechanical pressure and $S^{\circ}:= S + P\, I$, with $I$ denoting
the identity tensor, is the traceless part of $S$. We will assume throughout this paper that the material in all phases
does not support local densities for angular momentum (so-called couples). Hence the balance for angular momentum implies
that all stress tensors which appear are symmetric; note that all stress tensors are formulated in the embedding three-dimensional Euclidean space, i.e.\ are symmetric $3\times 3$-tensors.

Evidently, the simplest closure for the entropy flux in order to fulfill the entropy principle is $\Phi := \frac q T$,
which is the standard choice for single component materials. This leads to the reduced entropy production, being the
desired sum of binary products. Exploiting the symmetry of $S$, we obtain
\begin{align}
\label{entropy-prod-bulk}
\zeta = - \frac{1}{T}(\rho e + P -\rho s T -  \rho \mu)\, \div v
+ q \cdot \nabla \frac{1}{T}
 + \frac{1}{T} S^{\circ}:D^{\circ},
\end{align}
where $D:=\frac 1 2 (\nabla v + (\nabla v)^{\sf T} )$ is the symmetric part of the velocity gradient and
$D^{\circ}$ its traceless part. The dissipative mechanisms associated with these binary products are
``volume variations'', ``heat conduction'' and ``viscous shear'', in the order of their appearance in \eqref{entropy-prod-bulk}.
The simplest linear (in the co-factor) closure without cross-effects leads to the relations
\begin{align}
& \rho e + P -\rho s T -  \rho \mu = - \lambda \, \div v & \mbox{with } \lambda \geq 0,\label{P-closure}\\[1ex]
& q = \alpha \nabla \frac{1}{T}  & \mbox{with } \alpha  \geq 0,\\[1ex]
& S^{\circ} = 2 \eta D^{\circ} & \mbox{with } \eta \geq 0.\label{S-bulk}
\end{align}
Note that the closure parameters $\lambda, \alpha, \eta$
are allowed to depend on the basic variables, say $(\rho, T)$.
Hence, in particular, the heat flux closure is equivalent to Fourier's law. For consistency with standard notation,
we use $2 \eta$ instead of $\eta$, above. At this point, an explanation concerning \eqref{P-closure} is at order:
The only quantity which requires a closure is $P=-\frac 1 3 {\rm tr}\, S$. For a stagnant fluid, equation \eqref{P-closure}
reduces to $\rho e + p -\rho s T -  \rho \mu = 0$, where $p$ denotes the pressure at equilibrium.
We therefore let the thermodynamic pressure $p$ be defined by the Gibbs-Duhem relation, i.e.\ by
\begin{equation}
\label{GD-bulk}
\rho e + p -\rho s T =  \rho \mu .
\end{equation}
Then $P=p+\pi$ with the non-equilibrium pressure contribution $\pi$ and \eqref{P-closure} becomes
\begin{equation}
\label{pi-bulk}
\pi= - \lambda \, \div v.
\end{equation}
The irreversible pressure contribution $\pi$ is due to volume variations and the linear closure to model it reads as
$\pi = - \lambda \, \div v$. Let us note that the thermodynamic pressure $p$ from \eqref{GD-bulk} satisfies the Maxwell relation
$p=\rho^2 \frac{\partial \psi}{\partial \rho}$ with the free energy $\psi=\psi(T,\rho) :=e-s T$. Alternatively, one can define $p$
by the latter relation and obtain the Gibbs-Duhem relation \eqref{GD-bulk}
as a consequence. Note also that the entropy production
\eqref{entropy-prod-bulk} can now be written more concisely as
\begin{equation}
\zeta = q \cdot \nabla \frac{1}{T}  + \frac{1}{T} S^{\rm irr}:D,
\end{equation}
where the irreversible stress part is defined as $S^{\rm irr} = -\pi I + S^\circ$, but it is important to notice that
the last term represents two independent binary products.
$\mbox{ }$\\[1ex]
{\bf interface.} The same line of arguments leads to
\begin{align}
\label{surface-GD}
\Phi^\Sigma = \frac{q^\Sigma}{T^\Sigma}\quad \mbox{ and } \quad
 \rho^\Sigma e^\Sigma + p^\Sigma -\rho^\Sigma s^\Sigma T^\Sigma = \rho^\Sigma \mu^\Sigma
\end{align}
as well as
\begin{align}
\label{surf-entropy-prod}
& \zeta^\Sigma  = q^\Sigma \cdot \nabla_\Sigma \frac{1}{T^\Sigma}
- \frac{1}{T^\Sigma} \pi^\Sigma \, \div_\Sigma v^\Sigma
 + \frac{1}{T^\Sigma} S^{\Sigma,\circ}: D^{\Sigma,\circ}\nonumber \\[1ex]
& + \frac{1}{T^\Sigma} [\![  (v-v^\Sigma)_{||} \cdot (S \cdot n_\Sigma)_{||} ]\!]
 + [\![ \Big( \frac{1}{T}- \frac{1}{T^\Sigma}\Big) \Big( \dot{m} (e + \frac{p}{\rho}) +q\cdot n_\Sigma \Big) ]\!] \\[1ex]
& - [\![ \Big(\frac{\mu}{T}- \frac{\mu^\Sigma}{T^\Sigma}
+ \frac{1}{T^\Sigma} \big(\frac{(v-v^\Sigma)^2}{2}-n_\Sigma \cdot \frac{S^{\rm irr}}{\rho} \cdot n_\Sigma \big)\Big) \dot{m} ]\!]; \nonumber
\end{align}
recall that $\dot{m}=\rho (v-v^\Sigma) \cdot n_\Sigma$.
Here $\pi^\Sigma$ is the irreversible part of the interface pressure defined via
$\pi^\Sigma + p^\Sigma = -\frac 1 2 {\rm tr}\, S^\Sigma$ with the thermodynamic interface pressure $p^\Sigma$ from (\ref{surface-GD})$_2$. Moreover, $D^\Sigma = \frac 1 2 I_\Sigma \big(\nabla_\Sigma v^\Sigma
+ (\nabla_\Sigma v^\Sigma)^{\sf T} \big) I_\Sigma $ is the symmetric interface velocity gradient,
$D^{\Sigma,\circ}$ its traceless part and $I_\Sigma = I - n_\Sigma \otimes n_\Sigma $ denotes the surface projector, also called surface identity.

The dissipative processes associated with the binary products in \eqref{surf-entropy-prod} are, in the order of their appearance,
interfacial heat conduction, area variation, interfacial shear, one-sided slip between the interface and a bulk phase,
heat transfer to and from the interface and, finally, mass transfer to and from the interface.
The following closure relations result by assuming linear relations between the corresponding co-factors with one exception:
the mass transfer to or from the interface, i.e.\ the ad- and desorption processes
$\dot{m} = \dot{m}^{\rm ad} - \dot{m}^{\rm de}$,
will be modeled using a non-linear
relationship in analogy to the modeling of chemical reactions; cf.\ \cite{BotheDreyer}.
\begin{align}
& q^\Sigma = \alpha^\Sigma \nabla_\Sigma \frac{1}{T^\Sigma} & \mbox{with } \alpha^\Sigma \geq 0,\label{heat-int}\\[1ex]
& \pi^\Sigma = - \lambda^\Sigma \, \div_\Sigma v^\Sigma & \mbox{with } \lambda^\Sigma \geq 0,\label{visc1-int}\\[1ex]
& S^{\Sigma, \circ} = 2 \eta^\Sigma D^{\Sigma, \circ} & \mbox{with } \eta^\Sigma \geq 0,\label{visc2-int}\\[1ex]
& \beta^\Sigma (v-v^\Sigma )_{||} + ( S n_\Sigma)_{||} = 0 & \mbox{with } \beta^\Sigma \geq 0,\label{slip-int}\\[1ex]
& \frac{1}{T}- \frac{1}{T^\Sigma} +\delta^\Sigma \big( \rho (e + \frac{p}{\rho})(v-v^\Sigma) +q \big) \cdot n_\Sigma =0 & \mbox{with } \delta^\Sigma \geq 0,\label{heat-trans}\\[1ex]
& a^\Sigma \ln  \frac{\dot{m}^{\rm ad}}{ \dot{m}^{\rm de}}  =
\frac{\mu}{T}- \frac{\mu^\Sigma}{T^\Sigma}
+ \frac{1}{T^\Sigma} \big(\frac{(v-v^\Sigma)^2}{2}-n_\Sigma \cdot \frac{S^{\rm irr}}{\rho} \cdot n_\Sigma \big)
& \mbox{with } a^\Sigma \geq 0.\label{sorption-int}
\end{align}
The closure relation \eqref{sorption-int} employs the decomposition $\dot m = \dot{m}^{\rm ad}-\dot{m}^{\rm de}$.
Note that  \eqref{sorption-int} only fixes the ratio of ad- and desorption, while one of the rates
needs to be modeled based on experimental knowledge or a micro-theory. The simplest choice is to assume a desorption
rate according to $ \dot{m}^{\rm de} = k^{\rm de} \rho^\Sigma$ with $k^{\rm de} >0$.
Observe also that \eqref{sorption-int} is an implicit
equation regarding $ \dot{m}$, since $v^\pm -v^\Sigma=(v^\pm -v^\Sigma)_{||} +  \dot{m}^\pm / \rho^\pm$.

At this point it should be noted that the closure relations above are given in a condensed notation:
relations \eqref{heat-int}, \eqref{visc1-int} and \eqref{visc2-int} are employed for every interface $\Sigma^k$ ($k=1, 2, 3$)
separately with respective transport coefficients, while the transmission relations \eqref{slip-int}, \eqref{heat-trans}
and \eqref{sorption-int} apply to each interface in combination with any of the two adjacent bulk phases. In total, the closure
relations hence yield nine conditions at each of the three interfaces.
\\[1ex]
{\bf triple line.} Since the triple line is one-dimensional, the contact line stress tensor satisfies
$S^{\mathcal C}=-P^{\mathcal C} I_{\mathcal C}$ with the mechanical line pressure $P^{\mathcal C}:=-{\rm tr}\, S^{\mathcal C}$
and the line projector defined by $I_{\mathcal C} w =\langle w, \tau \rangle \tau$ with $\tau$ a unit tangent field
on ${\mathcal C}$.
By the same procedure as above, we obtain the following identities, where
$P^{\mathcal C}=p^{\mathcal C}+\pi^{\mathcal C}$, and $\dot{m}^\Sigma = \rho^\Sigma(v^\Sigma-v^{\cC})\cdot N$,
i.e.\ $\dot{m}^\Sigma_k = \rho^\Sigma_k (v^\Sigma_k -v^{\mathcal C})\cdot N^k$ for $k=1,2,3$.
\begin{align}
\Phi^{\mathcal C} = \frac{q^{\mathcal C}}{T^{\mathcal C}} \quad \mbox{ and } \quad
\rho^{\mathcal C} e^{\mathcal C} + p^{\mathcal C} -\rho^{\mathcal C} s^{\mathcal C} T^{\mathcal C} =  \rho^{\mathcal C} \mu^{\mathcal C}
\end{align}
as well as
\begin{align}
& \zeta^{\mathcal C} = q^{\mathcal C} \cdot \nabla_{\mathcal C} \frac{1}{T^{\mathcal C}}
 - \frac{1}{T^{\mathcal C}} \pi^{\mathcal C} \, \div_{\mathcal C} v^{\mathcal C}
+ \frac{1}{T^{\mathcal C}} [\![\![  (v^\Sigma-v^{\mathcal C})_{|||} \cdot (S^\Sigma \cdot N)_{|||} ]\!]\!] \nonumber\\[1ex]
& + [\![\![ \Big( \frac{1}{T^\Sigma}- \frac{1}{T^{\mathcal C}}\Big)
\Big( \rho^\Sigma (e^\Sigma + \frac{p^\Sigma}{\rho^\Sigma})(v^\Sigma-v^{\mathcal C}) +q^\Sigma \Big) \cdot N ]\!]\!] \\[1ex]
& - [\![\![ \Big(\frac{\mu^\Sigma}{T^\Sigma}- \frac{\mu^{\mathcal C}}{T^{\mathcal C}}
+ \frac{1}{T^{\mathcal C}} \big(\frac{(v^\Sigma-v^{\mathcal C})^2}{2}-N \cdot \frac{S^{\Sigma, \rm irr}}{\rho^\Sigma} \cdot N \big)\Big)\dot{m}^\Sigma ]\!]\!]. \nonumber
\end{align}
Above, the notation $(\cdot )_{|||}$ denotes the component tangential to the triple line and
$S^{\Sigma, \rm irr}:= - \pi^\Sigma I_\Sigma + S^{\Sigma, \circ}$.
In analogy to the interface we obtain the following closure relations for the dissipative processes on the triple line.
\begin{align}
& q^{\mathcal C} = \alpha^{\mathcal C} \nabla_{\mathcal C} \frac{1}{T^{\mathcal C}} & \mbox{with } \alpha^{\mathcal C} \geq 0,\label{heat-CL}\\[1ex]
& \pi^{\mathcal C} = - \lambda^{\mathcal C} \, \div_{\mathcal C} v^{\mathcal C} & \mbox{with } \lambda^{\mathcal C} \geq 0,\label{visc-CL}\\[1ex]
& \beta^{\mathcal C} (v^\Sigma-v^{\mathcal C} )_{|||} + ( S^\Sigma N)_{|||} = 0 & \mbox{with } \beta^{\mathcal C} \geq 0,\label{slip-CL}\\[1ex]
& \frac{1}{T^\Sigma}- \frac{1}{T^{\mathcal C}} +\delta^{\mathcal C} \big( \rho^\Sigma (e^\Sigma + \frac{p^\Sigma}{\rho^\Sigma})(v^\Sigma-v^{\mathcal C}) +q^\Sigma \big) \cdot N =0& \mbox{with } \delta^{\mathcal C} \geq 0,\label{trans-CL}\\[1ex]
& a^{\mathcal C} \ln  \frac{\dot{m}^{\Sigma, \rm ad}}{ \dot{m}^{\Sigma, \rm de}}  =
\frac{\mu^\Sigma}{T^\Sigma}- \frac{\mu^{\mathcal C}}{T^{\mathcal C}}
+ \frac{1}{T^{\mathcal C}} \big(\frac{(v^\Sigma-v^{\mathcal C})^2}{2}-N \cdot \frac{S^{\Sigma, \rm irr}}{\rho} \cdot N \big)
& \mbox{with } a^{\mathcal C} \geq 0.\label{sorption-CL}
\end{align}
As in the interface case, in \eqref{sorption-CL} the decomposition of $\dot{m}^{\Sigma}=\rho^\Sigma (v^\Sigma -v^{\mathcal C})\cdot N$ as $\dot{m}^{\Sigma}=\dot{m}^{\Sigma, \rm ad}- \dot{m}^{\Sigma, \rm de}$ is employed.
The relation \eqref{sorption-CL} governs the ratio of ad- and desorption at the triple line,
while one of the rates needs to be modeled based on experimental knowledge or a micro-theory.
Below we assume the desorption rate to be given by
$ \dot{m}^{\Sigma, \rm de} = k^{\Sigma, \rm de} \rho^{\mathcal C}$ with $k^{\Sigma, \rm de} >0$.
The transfer relations \eqref{slip-CL}, \eqref{trans-CL} and \eqref{sorption-CL} exist for every combination of the
triple line with one of the interfaces, of course with individual transfer coefficients.

To complete the model it remains to fix free energy functions for the bulk phases, the interfaces and the triple line.
This will only be done for a reduced model below.

\section{Isothermal Case with Vanishing Triple Line Mass}
\label{isothermal-model}
We consider the limiting case of isothermal conditions, i.e.\ the internal energy balances are replaced by a known constant temperature field; in particular, we have $T_{|\Sigma}=T^\Sigma$ and $T^\Sigma_{|{\mathcal C}}=T^{\mathcal C}$.
We also reduce the model complexity by neglecting the mass and inertia on the triple line. Moreover, we neglect any irreversible
stress contributions both on the interfaces and on the triple line. For consistency with the notation in interfacial
science, we do not employ the surface and line pressure, but rather let
$S^\Sigma = \gamma^\Sigma I_\Sigma$ and $S^{\mathcal C} = \gamma^{\mathcal C} I_{\mathcal C}$ with the interface tensions $\gamma^\Sigma$ and the line tension $\gamma^{\mathcal C}$.
Observe that this means $S^{\Sigma, irr}=0$. Because of zero triple line mass and isothermal conditions,
we assume the line tension to be constant.\\[1ex]
{\bf bulk phase.}
\begin{align}
& \partial_t \rho + \div (\rho v) =0,  \label{eq1}\\[0.5ex]
& \partial_t (\rho v) + \div (\rho v\otimes v) = \div S + \rho b, \label{eq2}\hspace{3.2in}
\end{align}
where the stress is given by $S= (-p + \lambda \div v )I + 2 \eta D^{\circ}$ according to \eqref{S-bulk} and \eqref{pi-bulk}.
In the compressible case, an equation of state
in the form $p=p(\rho)$ (with a strictly increasing function $p(\cdot )$)
is to be added according to the specific fluid under consideration.\\[1ex]
{\bf interface.} We again use the abbreviation $\dot{m}=\rho(v-v^\Sigma)\cdot n_\Sigma$. Then
\begin{align}
& \partial_t^\Sigma \rho^\Sigma + \div_\Sigma (\rho^\Sigma v^\Sigma) +[\![ \dot{m} ]\!] =0, \label{eq3}\\[1ex]
& \partial_t^\Sigma (\rho^\Sigma v^\Sigma)+ \div_\Sigma (\rho^\Sigma v^\Sigma\otimes v^\Sigma )
+[\![  v \, \dot{m}  ]\!] =[\![  S\cdot n_\Sigma ]\!]  + \div_\Sigma S^\Sigma + \rho^\Sigma b^\Sigma. \label{eq4}\hspace{0.4in}
\end{align}
Note that
\begin{align}
\div_\Sigma S^\Sigma  = \gamma^\Sigma \kappa_\Sigma n_\Sigma
+\nabla_\Sigma \gamma^\Sigma
\end{align}
in the considered case without surface viscosities.

This is complemented by the constitutive transmission conditions
\begin{align}
& \beta^\Sigma (v-v^\Sigma )_{||} + ( S n_\Sigma)_{||} = 0,\label{slip-law}\\[1ex]
& a^\Sigma \ln  \frac{\dot{m}^{\rm ad}}{ \dot{m}^{\rm de}}  =
\mu -  \mu^\Sigma + \frac{(v-v^\Sigma)^2}{2}- n_\Sigma\cdot \frac{S^{irr}}{\rho}\cdot n_\Sigma, \label{sorp-law}\hspace{2.2in}
\end{align}
where $a^\Sigma, \beta^\Sigma \geq 0$. In addition, the material dependent interface free energy function is required.
The latter determines
especially the interfacial equation of state $\gamma^\Sigma=\gamma^\Sigma(\rho^\Sigma)$ and we assume that $\gamma^\Sigma$
is a strictly decreasing function (i.e., the interface pressure depends strictly increasing on the interface mass density).
\\[1ex]
{\bf triple line.} In analogy with the interface-related notation, we use as before the abbreviation
\[
\dot{m}^\Sigma := \rho^\Sigma (v^\Sigma -v^{\mathcal C})\cdot N,
\]
i.e.\ $\dot{m}^\Sigma_k = \rho^\Sigma_k (v^\Sigma_k -v^{\mathcal C})\cdot N^k$ for $k=1,2,3$.
Due to $\rho^{\mathcal C} \equiv 0$, the triple line mass and momentum balances become
\begin{align}
& [\![\![ \dot{m}^\Sigma ]\!]\!] =0, \label{CL-Kirchhoff}\\[1ex]
& [\![\![ v^\Sigma \dot{m}^\Sigma ]\!]\!]=  [\![\![ \gamma^\Sigma N]\!]\!] + \gamma^{\mathcal C} \div_{\mathcal C} I_{\mathcal C}.
\label{CL-momKirchhoff}
\hspace{3.3in}
\end{align}
This is complemented by the constitutive transmission conditions
\begin{align}
& v^{\Sigma}_{1,|||} = v^{\Sigma}_{2,|||} =v^{\Sigma}_{3,|||}=: v_{|||}^\cC ,\label{CL-transfer}\\
& 
 \mu^{\Sigma}_k  -\mu^{\cC}+  ((v_k^{\Sigma}-v^{\mathcal C})\cdot N^k)^2/2=0 \qquad (k=1,2,3).\label{CL-transcond}\noindent
\end{align}
A few comments are at order: for simplicity, we consider the no-slip
condition $(v^\Sigma-v^{\mathcal C})_{|||}=0$, but note that the barycentric triple line velocity $v^{\mathcal C}$ is undefined for a triple line with zero mass. We consider  $v^{\mathcal C}$ as the kinematic velocity of the contact line $\cC$.
The chemical potential $\mu^{\mathcal C}$ is determined by one of the equations in \eqref{CL-transcond},
so actually only two equations remain there.
Also, observe that $v^\Sigma-v^\cC \perp \tau,n_\Sigma$, hence
\begin{equation}\label{v-perp}
 v^\Sigma-v^\cC = ( v^\Sigma-v^\cC | N )\, N,
\end{equation}
as $\CC (t) \subset \Sigma_k (t)$ for all times.

\section{Thermodynamical Consistency and Equilibria}
For this reduced isothermal model we show that the {\em total available energy}, i.e.\ the sum of the total kinetic energy and the total free energy
is a strict Lyapunov function in case of vanishing body forces. We hence let
\begin{equation}
\label{E(t)}
{\sf E}_a(t)=\int_{G} \rho (\frac{v^2}{2} +\psi) \,dx
+ \int_{\Sigma} \rho^\Sigma (\frac{(v^\Sigma)^2}{2} +\psi^\Sigma) \,do
+ \int_{\mathcal{C}} \gamma^{\mathcal C} \,dl,
\end{equation}
where $G$ is the total domain. We are going to show
\begin{theorem}
Let $(\rho, v,\rho^\Sigma, v^\Sigma, \Sigma, \mathcal{C})$ be a classical solution of the model from Section~\ref{isothermal-model}, i.e.\ a classical solution to
\eqref{eq1}, \eqref{eq2} with $S= (-p + \lambda \, \div v )I + 2 \eta D^{\circ}$, where $p(\rho)$ is strictly increasing in $\rho$,
$\lambda, \eta >0$ and $b=0$, \eqref{eq3}, \eqref{eq4} with $S^\Sigma = \gamma^\Sigma I_\Sigma$, where $\gamma^\Sigma (\rho^\Sigma)>0$ is strictly decreasing in $\rho^\Sigma$, and $b^\Sigma =0$,
\eqref{slip-law} with $\beta^\Sigma >0$, \eqref{sorp-law} with $a^\Sigma >0$,
\eqref{CL-Kirchhoff}, \eqref{CL-momKirchhoff} with $\gamma^{\mathcal C}$ a positive constant,
\eqref{CL-transfer} and \eqref{CL-transcond}.
We also assume that this solution is non-degenerate at the contact line, i.e.\ the interfaces meet at angles different from $0$ or $\pi$.
At the outer boundary, we assume $v\cdot n=0$, $v\cdot S n= 0$ on $\partial G$ and $v^\Sigma \cdot N=0$ on $\Sigma \cap \partial G$.\\
Then the total available energy ${\sf E}_a$ from \eqref{E(t)} is a strict Lyapunov function.
\end{theorem}
{\bf Proof.} (i) Let $(\rho, v,\rho^\Sigma, v^\Sigma, \Sigma, \mathcal{C})$ be a classical solution of the model from Section~\ref{isothermal-model}.
For the bulk contribution, we first apply the transport relation \eqref{vol-trans2} and use the momentum balance \eqref{bulk-momentum}
to eliminate $\rho \frac{Dv}{Dt}$. We then exploit $\psi=\psi(\rho)$ with $\psi' (\rho)=p / \rho^2$ and use the mass
balance \eqref{bulk-mass} to eliminate $\frac{D\rho}{Dt}$. Application of the two-phase divergence theorem for partial
integration in the form
\[
\int_G v \cdot div S \, dx = \int_{\partial G} v\cdot Sn \, do - \int_G S:\nabla v \, dx
-\int_\Sigma [\![ v\cdot Sn_\Sigma ]\!] \, do
\]
yields
\begin{align}
\frac{d}{dt} \int_{G} \rho (\frac{v^2}{2} +\psi) \,dx & =
- \int_{\partial G} \rho (\frac{v^2}{2} +\psi) v \cdot n \, do + \int_{\partial G} v \cdot Sn \, do \nonumber \\[1ex]
& -\int_G S^{\rm irr} : \nabla v \, dx
-\int_\Sigma [\![ v \cdot S n_\Sigma ]\!] \, do
+\int_\Sigma [\![ \dot m (\psi + \frac{v^2}{2} ) ]\!] \, do.\label{id-bulk}
\end{align}
For the interface contribution, we first apply the transport relation \eqref{surf-trans2} and use the momentum balance
in the non-conservative form
\[
\rho^\Sigma \frac{D^\Sigma v^\Sigma}{Dt} +
[\![ (v-v^\Sigma)\dot m]\!] = [\![ S \cdot n_\Sigma]\!] + \div_\Sigma S^\Sigma,
\]
which follows from \eqref{surf-momentum} and \eqref{surf-mass}, to eliminate $\rho^\Sigma \frac{D^\Sigma v^\Sigma}{Dt}$.
Next, we apply the surface divergence theorem for partial integration of $v^\Sigma \cdot \div_\Sigma S^\Sigma$,
employ $\psi^\Sigma= \psi^\Sigma (\rho^\Sigma)$ with
$(\psi^\Sigma)'(\rho^\Sigma)=p^\Sigma / (\rho^\Sigma)^2=- \gamma^\Sigma / (\rho^\Sigma)^2$, \eqref{surf-mass} and the interface
Gibbs-Duhem relation (\ref{surface-GD})$_2$ to obtain
\begin{align}
& \frac{d}{dt}  \int_{\Sigma} \rho^\Sigma (\frac{(v^\Sigma)^2}{2} +\psi^\Sigma) \,do =
\int_{\partial \Sigma} \rho^\Sigma  (\frac{(v^\Sigma)^2}{2}+ \psi^\Sigma) (V_{\partial \Sigma} - v^\Sigma\cdot N )\,dl
+ \int_{\partial \Sigma} v^\Sigma  S^\Sigma \cdot N \,dl\nonumber \\[1ex]
& - \int_\Sigma S^{\Sigma, \rm irr} : D^\Sigma \,do
+ \int_\Sigma v^\Sigma \cdot  [\![ S n_\Sigma ]\!]  \,do
+ \int_\Sigma [\![ (\frac{(v^\Sigma)^2}{2} - v^\Sigma \cdot v -\mu^\Sigma) \dot m ]\!]\, do.\label{id-surf}
\end{align}
For the triple line contribution, we apply the transport relation \eqref{CL-trans}. For constant $\gamma^{\mathcal{C}}$ this yields
\begin{align}
& \frac{d}{dt}  \int_{\mathcal{C}} \gamma^{\mathcal C} \,dl =
\int_{\mathcal{C}} \gamma^{\mathcal C} \div_{\mathcal C} v^{\mathcal C} \,dl =
\int_{\mathcal{C}} \gamma^{\mathcal C} I_{\mathcal C} : \nabla_{\mathcal C} v^{\mathcal C} \,dl =
- \gamma^{\mathcal C} \int_{\mathcal{C}} v^{\mathcal C} \cdot \div_{\mathcal C} I_{\mathcal C} \,dl.\label{id-CL}
\end{align}
Employing the identi\-ties \eqref{id-bulk}, \eqref{id-surf} and \eqref{id-CL}, we obtain
\begin{align}
\label{total-free-energy1}
\dot{{\sf E}}_a = & \int_{\partial \Sigma} \rho^\Sigma  (\frac{(v^\Sigma)^2}{2}+ \psi^\Sigma) (V_{\partial \Sigma} - v^\Sigma\cdot N )\,dl
+ \int_{\partial \Sigma} v^\Sigma  S^\Sigma \cdot N \,dl\nonumber \\[1ex]
& -\int_G S^{\rm irr} : D \, dx - \int_\Sigma S^{\Sigma, \rm irr} : D^\Sigma \,do
-\int_\Sigma [\![ (v - v^\Sigma)_{||} \cdot (S^{\rm irr} n_\Sigma)_{||} ]\!] \, do\nonumber \\[1ex]
& + \int_\Sigma [\![ (\mu -\mu^\Sigma + \frac{(v-v^\Sigma)^2}{2} - n_\Sigma \cdot \frac{S^{\rm irr}}{\rho}\cdot n_\Sigma)\, \dot m ]\!]\, do
 - \gamma^{\mathcal C} \int_{\mathcal{C}} v^{\mathcal C} \cdot \div_{\mathcal C} I_{\mathcal C} \,dl.
\end{align}
Inserting the constitutive relation $S^\Sigma = \gamma^\Sigma I_\Sigma$ and exploiting the assumptions
$v\cdot n=0$, $v\cdot S n= 0$ on $\partial G$ and $v^\Sigma \cdot N=0$ on $\Sigma \cap \partial G$, we get
\begin{align}
\label{total-free-energy2}
\dot{{\sf E}}_a = & -\int_G S^{\rm irr} : D \, dx -\int_\Sigma [\![ (v - v^\Sigma)_{||} \cdot (S^{\rm irr} n_\Sigma)_{||} ]\!] \, do\nonumber \\[1ex]
& + \int_\Sigma [\![ (\mu -\mu^\Sigma + \frac{(v-v^\Sigma)^2}{2} - n_\Sigma \cdot \frac{S^{\rm irr}}{\rho}\cdot n_\Sigma)\, \dot m ]\!]\, do\nonumber \\[1ex]
& + \int_{\mathcal{C}} [\![\![ (\frac{(v^\Sigma)^2}{2}+ \psi^\Sigma) \dot{m}^\Sigma ]\!]\!]\,dl
- \int_{\mathcal{C}} [\![\![  \gamma^\Sigma v^\Sigma \cdot N ]\!]\!]\,dl
- \gamma^{\mathcal C} \int_{\mathcal{C}} v^{\mathcal C} \cdot \div_{\mathcal C} I_{\mathcal C} \,dl.
\end{align}
Expanding $[\![\![  \gamma^\Sigma v^\Sigma \cdot N ]\!]\!]$ as $[\![\![  \gamma^\Sigma (v^\Sigma - v^{\mathcal C} + v^{\mathcal C}) \cdot N ]\!]\!]$, exploitation of \eqref{CL-momKirchhoff} allows to rewrite the triple line contribution as
\[
\int_{\mathcal{C}} [\![\![  (\frac{(v^\Sigma)^2}{2}+ \mu^\Sigma) \dot{m}^\Sigma ]\!]\!]\,dl
-  \int_{\mathcal{C}} v^{\mathcal C} \cdot [\![\![  v^\Sigma \dot{m}^\Sigma ]\!]\!] \,dl.
\]
Using $ (v^{\mathcal C})^2 [\![\![ \dot{m}^\Sigma ]\!]\!] =0$ due to \eqref{CL-Kirchhoff}, where
$v^{\mathcal C}_{|||}$ is given as the well-defined tangential part $v^{\Sigma}_{|||}$ by \eqref{CL-transfer},
we see that \eqref{total-free-energy2} implies
\begin{align}
\label{total-free-energy3}
\dot{{\sf E}}_a = & -\int_G S^{\rm irr} : D \, dx -\int_\Sigma [\![ (v - v^\Sigma)_{||} \cdot (S^{\rm irr} n_\Sigma)_{||} ]\!] \, do\nonumber \\[1ex]
& + \int_\Sigma [\![ (\mu -\mu^\Sigma + \frac{(v-v^\Sigma)^2}{2} - n_\Sigma \cdot \frac{S^{\rm irr}}{\rho}\cdot n_\Sigma)\, \dot m ]\!]\, do\nonumber \\[1ex]
& + \int_{\mathcal{C}} [\![\![ (\mu^\Sigma + \frac{(v^\Sigma - v^{\mathcal C})^2}{2}) \dot{m}^\Sigma ]\!]\!]\,dl.
\end{align}
To come to the final representation of $\dot{{\sf E}}_a$, we have to write out the jump brackets $[\![ \cdot ]\!]$
and $[\![\![ \cdot ]\!]\!]$. We start with the triple line contribution and have, by \eqref{CL-Kirchhoff},
\[
[\![\![ (\mu^\Sigma + \frac{(v^\Sigma - v^{\mathcal C})^2}{2}) \dot{m}^\Sigma ]\!]\!]
= - \sum_{k=1}^3 \big( \mu^\Sigma_k - \mu^{\mathcal C}
+ \frac{(v^\Sigma_k - v^{\mathcal C})^2}{2} \big) \dot{m}^\Sigma_k.
\]
Since $(v^\Sigma_k - v^{\mathcal C})^2=((v^\Sigma_k - v^{\mathcal C})\cdot N^k)^2$ by \eqref{v-perp}, the triple line contribution
vanishes due to the constitutive assumption \eqref{CL-transcond}.
Insertion of the other constitutive relations, i.e.\ \eqref{S-bulk}, \eqref{pi-bulk}, \eqref{slip-law} and \eqref{sorp-law},
finally leads to
\begin{align}
\label{total-free-energy4}
\dot{{\sf E}}_a = & -\int_G \lambda (\div v)^2 \, dx -\int_G 2 \eta D^\circ : D^\circ\, dx \nonumber \\[1ex]
& -\int_\Sigma \beta^{\Sigma, +} (v^+ - v^\Sigma)_{||}^2  \, do -\int_\Sigma \beta^{\Sigma, -} (v^- - v^\Sigma)_{||}^2  \, do\nonumber \\[1ex]
& - \int_\Sigma a^{\Sigma, +} \big( \log \dot{m}^{+,\rm ad} - \log \dot{m}^{+,\rm de}  \big)
\big( \dot{m}^{+,\rm ad} - \dot{m}^{+,\rm de} \big) \, do\nonumber \\[1ex]
& - \int_\Sigma a^{\Sigma, -} \big( \log \dot{m}^{-,\rm ad} - \log \dot{m}^{-,\rm de}  \big)
\big( \dot{m}^{-,\rm ad} - \dot{m}^{-,\rm de} \big) \, do.
\end{align}
Notice that, according to our condensed notation, the integrals over $\Sigma$ are to be taken over the three interfaces
and the notation $(\cdot)^\pm$ then denotes the respective one-sided bulk limits.
Evidently, \eqref{total-free-energy4} shows that ${\sf E}_a$ is decreasing along classical solutions, i.e.\ ${\sf E}_a$ is a Lyapunov function.


Next we want to characterize the equilibria of the problem, proving at the same time that the total available energy ${\sf E}_a$ is a strict Lyapunov functional for the system. To this end assume that we have a solution where ${\sf E}_a$ is not strictly decreasing at all times. Then there is an interval $J=(t_1,t_2)$ where ${\sf E}_a$ is constant, hence $d{\sf E}_a/dt=0$ in $J$. This implies, by \eqref{total-free-energy4},
$$ {\rm div}\,v=0,\quad D^\circ=0,\quad v^+_{||}=v^\Sigma_{||}=v^-_{||},\quad \dot{m}^+=\dot{m}^-=0,$$
as $\lambda, \eta, \beta^{\Sigma,\pm},a^{\Sigma,\pm}>0$ by assumption.
This yields $D=0$, as well as $[\![v]\!]=0$ on $\Sigma$, which by Lemma 1.2.1 of the monograph \cite{PrSi15} implies $v=v^{\Sigma}_{||}=0$. Next, investigating the equations for the bulk, we see that $\partial_t\rho=0$ and $\nabla p=0$,
which implies that $\rho$ is constant in the phases, as $p_k$ is by assumption a strictly increasing function of $\rho_k$.

In the next step, we look at the equations on the interfaces. By the definition of $\dot{m}^\pm$ we obtain
$$0 =\dot{m}^\pm = \rho^{\pm}( v^\pm-v^\Sigma)\cdot n^\Sigma = -\rho^{\pm} v^{\Sigma}\cdot n^\Sigma,$$
hence $  v^{\Sigma}\cdot n^\Sigma=0$ which yields $v^\Sigma=0$. Then the mass balance on $\Sigma$ implies $\partial_t^\Sigma \rho_\Sigma=0$ on $J$.
Furthermore, $v=0$ and $\rho$ constant yield $\mu^\pm$ constant, hence $\mu^\pm=\mu^\Sigma$ is constant by \eqref{sorp-law}. This shows that $\rho^\Sigma$ is constant, as $\mu^\Sigma$ is strictly increasing with $\rho^\Sigma$. To see the latter, recall that
$\mu^\Sigma = \psi^\Sigma + p^\Sigma / \rho^\Sigma$ and, hence, $(\mu^\Sigma)' (\rho^\Sigma) = (p^\Sigma)' (\rho^\Sigma) / \rho^\Sigma >0$.
This shows further that $\gamma^\Sigma$ is constant. Looking at the stress transmission condition this further yields $\kappa^\Sigma$ constant on each of the surfaces $\Sigma_k$; more precisely we obtain $\kappa^\Sigma = [\![p]\!]$, i.e.\ the Young-Laplace law holds on each of the surfaces $\Sigma_k$.

In the final step, we consider the equations on the contact line. Here we have $v^\cC_{|||}=0$ by \eqref{CL-transfer}, as well as
$$v^\cC = (v^\cC|N^k)N^k, \quad k=1,2,3,$$
hence $v^\cC=0$ if ${\rm dim}\, {\rm span}\{N^k\}_{k=1}^3=2$, i.e.\ in the non-degenerate case which is assumed to hold.
This further yields $\mu^\cC$ constant, and there remains the Kirchhoff law
$$ \sum_{k=1}^3 \gamma^kN^k =\gamma^\cC \kappa^\cC,$$
where $\kappa^\cC= -{\rm div}_\cC I_\cC= \nabla_\cC\tau$ denotes the curvature vector of the contact line $\cC$.

So, $\dot{{\sf E}}_a =0$ on $(t_1,t_2)$ implies the following:
\begin{enumerate}
\item The densities are constant, and all velocities vanish.
\item The curvatures $\kappa^{\Sigma_k}$ of the hypersurfaces $\Sigma_k$ are constant.
\item $\sum_{k=1}^3 \gamma^k N^k = \gamma^\cC \kappa^\cC$, where the coefficients $\gamma^j$ are positive constants.
\end{enumerate}
But this implies that the classical solution coincides with an equilibrium of the system at any $t\in (t_1,t_2)$,
hence remains at a fixed equilibrium for all $t>t_1$.
Consequently, it holds that for any classical solution,  ${\sf E}_a$ is strictly decreasing outside of equilibria,
 i.e.\ ${\sf E}_a$ is actually a strict Lyapunov function.\\
$\mbox{ }$\hfill $\Box$\\[2ex]
To identify all possible equilibrium configurations is a purely geometrical problem. It appears to be a challenging problem and will not be analyzed any further, here.\\[2ex]
\noindent
{\bf Final Remarks.} 1. The proof that ${\sf E}_a$ is non-increasing along classical solutions requires all interfacial
and triple line conditions, in particular the condition \eqref{CL-transcond}.
This confirms that the original interface formation model of
Shikhmurzaev misses one contact line condition. The origin of this transmission condition is the fact that transfer of
mass, here from one interface across the contact line to another interface, is a dissipative process which requires
a closure relation. This is similar to the case of mass transfer
across a fluid interface: even without interfacial mass, a fluid interface carries interfacial energy and, in general,
entropy can be produced at the interface. In order to avoid entropy production, the condition which guarantees zero
interfacial entropy production has to be added, leading in the simplest case to continuity of the chemical bulk potentials.
The triple line analog is equation \eqref{CL-transcond} above.
In the more general case of non-trivial entropy production, thermodynamically consistent closure
leads to a condition like \eqref{sorption-CL}.\\[0.5ex]
\indent
2. The molecular kinetic theory of dynamic contact lines
supports a friction-like dissipation term at the contact line, modeled as being
proportional to the square of the contact line speed. If, instead of the non-linear closure \eqref{sorption-CL},
a linear relation is imposed, the rate of entropy production due to transfer of interfacial mass across the
contact line becomes proportional to the contact line speed squared. In an isothermal setting, this entropy production
is proportional to the dissipation of available energy. Hence, the so-called contact line friction can be identified
with the interfacial mass transfer dissipation mechanism.
\section*{Acknowledgement}
The first author (D.B.) is grateful for support by the DFG within the cluster of excellence
``Center of Smart Interfaces'', TU Darmstadt. The second author (J.P.) thanks the DFG for continuous support
in the framework of individual research projects.

\end{document}